\documentclass[12pt ]{article} 

 \usepackage{pdfsync}    
 \usepackage{amsmath}
 \usepackage{verbatim}
 \usepackage{amsthm}
\newtheorem{definition}{Definition}[section]

\newtheorem{theorem}{Theorem}[section]
\newtheorem{lemma}{Lemma}[section]
\newtheorem{corollary}{Corollary}[section]

\newtheorem{remark}{Remark}[section]
\newcommand{\eqnsection}{
   \renewcommand{\theequation}{\thesection.\arabic{equation}}
   \makeatletter 
   \csname @addtoreset\endcsname{equation}{section}
   \makeatother}

\def \be{\begin{equation}}
\def \ee{\end{equation}}
\def \bt{\begin{theorem}}
\def \et{\end{theorem}}
\def \bea{\begin{eqnarray}}
\def \eea{\end{eqnarray}}
\def \bas{\begin{eqnarray*}}
\def \eas{\end{eqnarray*}}
\def \bl{\begin{lemma}} 
\def \el{\end{lemma}}

\def \al{\alpha}

\def \Ga{\Gamma}

\def \ff{\infty}
\def \wh{\widehat}

\def \cd{\,\cdot\,}

\def \DD{{\cal D}}

\def \({\left(}
\def \){\right)}

\def \nn{\nonumber}
\def \Proof{\noindent{\bf Proof $\,$ }}

\def \bc{\begin{center} }
\def \ec{\end{center} }
\def \bs{\begin{slide} }
\def \es{\end{slide} }

\def\square{{\vcenter{\vbox{\hrule height.3pt
        \hbox{\vrule width.3pt height5pt \kern5pt
           \vrule width.3pt}
        \hrule height.3pt}}}}
\def\qed{{\hfill $\square$ \bigskip}}

\eqnsection
 
\begin{document}

\def\wh{\widehat}
\def\ol{\overline}

\title{ON CRITICAL POINTS FOR GAUSSIAN VECTORS WITH INFINITELY DIVISIBLE SQUARES}
\author{Hana Kogan}
\maketitle
\begin{abstract}
The infinite divisibility of the squared Gaussian Process with non-zero mean is dependent on the signs of the row sums of the inverse of its covariance matrix. We give an explicit formula describing this relationship. 
Let the Gaussian vector have a critical point $\al_0$ if its squares are infinitely divisible for all $\al \leq \al_0$ and are not infinitely divisible for all $\al > \al_0$. We give the upper bound for the critical point of a non-zero mean Gaussian vector.
\end{abstract}
\section{Introduction }

Let $G = ( \eta_1, \eta_2,....\eta_n)$ be an $n$-dimensional Gaussian vector. We say that $G$ has infinitely divisible squares (or $G^2$ is infinitely divisible) if for any $m \in N$
 \be \label {0}
 G^2 :=(\eta_1^2, \eta_2^2,....\eta_n^2)\stackrel{law}{=} \sum_{i=1}^m{Z_i}
 \ee
  where $Z_i$ are independent identically distributed $n$-dimensional random  vectors.

  \begin{itemize}
 \item[]	  Let $N$ be an $n\times n$ matrix such that $N_{i,j}=0$  for all $i\neq j; N_{i,i}=\pm1$ for all $i$. $N$ is called a signature matrix. 
 \item[]	 A non-singular matrix $A$ is called an $M$-matrix if $A_{i,j}^{-1}\geq 0$ for all $1<i,j<n$  and $A_{i,j}\leq 0$  for all $i\neq j$. Note that $A_{i,i}$ must be positive.  
\item[]	  A symmetric matrix $A$ is irreducible if it can not be written as a direct sum of square matrices. 
  \end{itemize}
  
 	  In what follows all covariance matrices considered are irreducible. The ${\mathbf 1}$ denotes the $n$-dimensional vector with all entries equal to $1$. 
   The following Theorem is due to Griffiths and Bapat, see also \cite[Theorem 13.2.1]{book}  and \cite[Theorem 1.1]{inf div}. Theorem \ref{0-mean div} completely characterizes zero-mean Gaussian processes with infinitely divisible squares.  
   \begin{theorem}\label{0-mean div} 
  Let $G = ( G_1, G_2, ..... G_n)$ be a mean zero Gaussian random variable with a strictly positive-definite covariance matrix $\Gamma = \{\Gamma_{i,j}\} = \{E(G_iG_j)\}$.   $G^2$ is infinitely divisible if and only if $N\Gamma^{-1}N$ is an M-matrix for some signature matrix $N$. 
  \end{theorem}
  We are interested in Gaussian vectors with infinitely divisible squares. We can therefore restrict our attention to Gaussian vectors with positive covariance.
  
  The necessary and sufficient conditions for infinite divisibility of Gaussian squares with uniform non-zero mean - $$(G+\al {\mathbf 1})^2 := ( (G_1+\al)^2, (G_2+\al)^2,...)$$ was posed and answered by Eisenbaum and Kaspi,  \cite{EK}.
  Their result was developed to include the non-uniform mean Gaussian squares by Marcus and Rosen.
     Let ${\mathbf{c}}=(c_1, c_2...)$ be any real $n$-dimensional vector. The question about  infinite divisibility of $(G + \al{\mathbf{c}})^2:= ((G_1+\al c_1)^2, (G_2+\al c_2)^2,...)$ for all $\al \in  {\mathbf R}$ is answered by   Theorem 1.3, \cite{inf div}:
  \begin{theorem}\label{var mean}
  Let $G$ be a mean zero Gaussian vector with strictly positive definite covariance  matrix $\Gamma$, $N(0,b)$ denote the Gaussian random variable with mean $0$ and variance $b$. Let ${\mathbf{c}}=(c_1, c_2...)$ and $C$ be an $n\times n$ matrix with $C_{i,i}= {c_i}, C_{i,j}=0$ for $i\neq j$. The following are equivalent:
  \begin{enumerate}
 \item$(G+c\alpha)^2$ is infinitely divisible for all $\alpha \in {\mathbf{R}}$.
 \item There exists some $b\in{\mathbf{R}}, b\neq 0$ such that  for $\eta=N(0, b)$, independent of $G$, \quad $(G_1+c_1\eta, G_2+c_2\eta,...G_n+c_n\eta, \eta)$ has infinitely divisible squares. In this case  $(G_1+c_1\eta, G_2+c_2\eta,...G_n+c_n\eta, \eta)$ has infinitely divisible squares for all $b\in {\mathbf{R}}$.
 \item $C\Gamma^{-1}C$ is an $M$-matrix with non-negative row sums.
 \end{enumerate}
  \end{theorem}

  \begin{corollary}\label{nz}
  Let $G$, $\Gamma$, $C$ be as above. When 1, 2 or 3 of Theorem \ref{var mean} hold,
  \begin{enumerate}
  
   \item
   $C_{i,i}\neq 0$, for all $i$.
\item $ \Gamma_{i,j}\neq 0,\quad \forall (i,j)$.
\end{enumerate}
 \end{corollary}

 The natural question that could be asked at this point (and which was posed by Marcus and Rosen, see \cite{inf div} and \cite{critical}) is about the existence and characterization of the vectors $G$ for which $(G+\al{\mathbf{c}})^2$ is infinitely divisible for some but NOT ALL $\alpha$. We shall concentrate on the case ${\mathbf{c}}={\mathbf{1}}$.
 To put this question formally we need

 \begin{definition}\label{cr}
 Let $G$ be a Gaussian vector. Let $\al_0$ be a  real number $0< \al_{0}< \ff$. 
   We say that $\alpha_0$ is a critical point of $G$   if 
 \begin{itemize}
 \item 
 $(G+\alpha {\mathbf 1})^2$  is infinitely divisible for all   $|\alpha|\le \alpha_{0}$,

 and  $(G+\alpha {\mathbf 1})^2$  is not infinitely divisible for any $ |\alpha|> \alpha_{0}$ 

  \item  when   
$(G+\alpha {\mathbf 1})^2$ is  infinitely divisible for all $ \alpha $  then  $ G$ 
  has no critical point.  
 \end{itemize}
 \end{definition}

    Marcus and Rosen have shown in \cite{critical} the existence of a critical point $0<\alpha_0 < \ff$,   for all two-dimensional Gaussian vectors with infinitely divisible squares.  The question about existence of a critical points for the vectors of higher dimension remains open.
    
  We   explore the connection between the properties of the inverse covariance matrix $\Gamma^{-1}$ and the infinite divisibility of $(G+\al)^2$ for different $\al$. We  obtain some upper bounds on a critical point for Gaussian $n$-dimensional vectors $G$ with $n>2$ and also the proof of Theorem \ref{var mean} ( part $3 \Rightarrow 1 $) which does not invoke the Markov process theory (see Section 4). In particular, the original proof of  \ref{var mean} uses the Second Ray-Knight Theorem to establish an isomorphism between $(G+c)^2$ and the sum of two independent infinitely divisible  random vectors. We obtain an elementary proof of the same result using the well known result due to Griffith and Bapat,  (see \cite{book}), Theorem \ref{main1} and Lemma \ref{q pos}.
  
  \medskip
 We employ the necessary and sufficient conditions for infinite divisibility of a random vector due to Feller;(see \cite{F}, Ch. XIII.4), \cite{B},  \cite[13.2.2]{book},  and the application therein to Gaussian vectors. We follow the presentation of Marcus and Rosen in \cite[Lemmas 13.3.1 and 5.2.1] {book}, and \cite{critical}:
\bl \label{div cond}\cite[Lemma 13.2.2]{book}
  Let $G$ be a Gaussian vector with positive definite covariance $\Ga$, $\al \in \mathbf{R}$.
For any vector $\tilde{\Lambda} \in {\mathbf{R}^n}$ denote by $\Psi_G(\tilde{\Lambda})$  the Laplace transform of $(G+\alpha)^2$ and let $S$ be a diagonal $n \times n$  matrix with entries $s_{i,i}\in(0,1]$. Let $t>0$.  Let $\Lambda$ be a diagonal matrix such that $\Lambda_{i,i}= \tilde{\Lambda}_i$ for all $i$. For any $\tilde{\Lambda} \in {\bf{R}^n}$ we can write
\be \label{4}
\Lambda=t(I-S)
\ee for any $t$ sufficiently large and some $S$ , where $S$ is a diagonal  $n\times n$ matrix with $S_{i,i} \in [0,1]$.
 Let $$\Phi (t, S)=\log{\Psi_G(\tilde{\Lambda})}$$ and suppose that $\Phi (t, S)$ has a power series expansion about $S = {\mathbf{0}}$.
Then $(G+\alpha)^2$ is infinitely divisible if and only if for all $t$ sufficiently large  all coefficients of this expansion are non-negative, except the constant term. 
\el  
  
 We obtain the series expansion for  $\Phi (t, S)$ for general $\al$  in a form that makes it relatively easy to analyze the sign of the coefficients.
 
  Let  \be\label{8} Q(t)=[I+(t\Gamma)^{-1}]^{-1} = I-(I+t\Gamma)^{-1}.\ee
  Let $D \in \mathbf{R}^n$ such that 
  \begin{equation} \label{1}
  D_i = \sum_{k=1}^n \Gamma^{-1}_{i,k} ; 
  \end{equation}
   i.e $D_i$ denotes the $i$-th row sum of $\Gamma^{-1}$.
   
  Let
$$\Phi(t,S)= \Phi_1(t,S)+\Phi_2(t,S);$$
where
  $$\Phi_1(t,S)= \Phi(t,S)|_{\al=0}$$
The following series expansion for $\Phi_1(t,S)$ is the result of Bapat and Griffith:
 
    \be \label{11}
    \Phi_1(t,S)=\frac{1}{2} \Big(\log |I-Q| + \sum^{\infty}_{m=1}\frac{{\rm trace }\,(QS)^m}{m }\Big). \ee
 We extend their result to include $\al \in R$:

 \begin{theorem}\label{main1}
  
Under the notation as above,

$$\Phi(t,S)= \Phi_1(t,S)+\Phi_2(t,S);$$
and for all $\al \in R$
 
 \be \label{18}   
  \Phi_2 (S)=\frac{\alpha^2}{2} \(-{\mathbf 1}\!\cdot\!\Ga^{-1}\!\cdot\! {\mathbf 1}^T + D\!\cd \!\( \frac{1}{t} \sum_{m=1}^{\ff}   (QS)^m  \)\!\cd \!D^T\)(1+O(1/t)).  
\ee

 \end{theorem}

Since $(G+\al)^2$ is infinitely divisible if and only if the series expansion for $\Phi(S)$ has non-negative coefficients for all non-constant terms,
it is clear from Theorem \ref{main1} that the signs of coefficients in the expansion of $\Phi_2(S)$- and hence the infinite divisibility of $(G+\al)^2$  will depend on the signs of the entries of $Q$ and $D$. In particular, since the coefficients of $\Phi_2$ depend linearly on $\al^2$, it is clear that if $(G+\al)^2$ is not infinitely divisible for a particular $\al$ then it is not infinitely divisible for all $\beta>\al$, and so the Definition \ref{cr} makes sense.
The following Lemma establishes the positivity of the entries of $Q$, which makes possible a relatively straightforward analysis of the coefficients of $\Phi(S)$. It also provides an alternative, elementary proof of Theorem \ref{var mean} which does not invoke the Isomorphism Theorems and Markov Processes.
\begin{lemma} 

 \label{q pos}
For all $t$ sufficiently large,
$$Q >0$$
\end{lemma}   
 
 Theorem \ref{main1} and Lemma \ref{q pos}  allow us to obtain an upper bound for the critical point of $G$:
  \begin{theorem}\label{my thm} 
  Let $G$ be a mean zero Gaussian vector with infinitely divisible squares and covariance matrix   $\Ga$ (hence $\Ga^{-1}$ is an $M$ matrix).
   
  Let
\be \label{2} \mathcal{D} = \{ (i,j) : D_iD_j <0\}.
\ee

If $\DD\ne \emptyset $   and $\al_{0}$  is a critical point of $G$,  then 
  \begin{equation}\label{3} 
 \alpha_0\leq \inf_{\mathcal{D}}{\left\{\(\frac{\Gamma^{-1}_{i,j}}{2D_iD_j}\)^{\frac{1}{2}}\right\}}. \end {equation} 
 \end{theorem}
 
 An immediate consequence of Theorem \ref{my thm} is
 
 \begin{corollary}\label{0 infty}

 \begin{itemize}
\flushleft
 \item[ $\bullet$] If for some pair $(i,j) \in \mathcal{D} \quad \Ga^{-1}_{i,j} =0,$ then $\al_0 =0$.
 \item [ $\bullet$]if $\DD= \emptyset $ then $G$ has no critical point.
  
 \end{itemize}
 
\end{corollary}

Note that for the $\mathcal{D}$ to be non-empty $\Ga^{-1}$ matrix has to have  at least one negative row sum. Hence the second result of the Corollary is the statement   (  $3 \Rightarrow 1 $) of Theorem \ref{var mean}. The first is the new result that allows us to construct examples of Gaussian vectors with infinitely divisible squares and zero critical point (see Section 7).
 We prove Theorem \ref{main1} and Lemma \ref{q pos} in Section 2. Theorem \ref{my thm} and Corollary \ref{0 infty} are proved in Section 3. Section 4 contains an elementary proof of Theorem \ref{var mean}. Various applications of Theorem \ref{main1} are given in Sections 5 and 6. Section 7 contains some numerical examples.

 \section{Proofs of Theorem \ref{main1} and Lemma \ref{q pos} . }
We use the following
\begin{lemma}\label{phi2}

 \be
\Phi_2 (S) =  \frac{\alpha^2 t}{2} 1\cdot[(Q-I)+ (I-Q^{-1})\sum_{m=1}^\infty{(QS)^m}(Q-I)]\cdot1^T\label{12}
\ee   
\end{lemma} 
 \proof
By \cite{book}, Lemma 5.2.1, \begin{equation}\label{5}
                                      \Psi_G(\tilde{\Lambda})=\det(I+\Gamma\Lambda)^{-1/2}\cdot(\exp(\frac{\alpha^2}{2}1[\Lambda\tilde{\Gamma}\Lambda-\Lambda]1^T)). \end{equation}
where $\tilde{\Gamma}=(\Gamma^{-1}+\Lambda)^{-1}$.

 Let 
 \be
 \quad \Lambda  =t(I-S)\nn
 \ee
  for sufficiently large $t$.

 Then \begin{equation}\label{7}
  \Phi_2(t,S)=\frac{\alpha^2}{2}1[\Lambda\tilde{\Gamma}\Lambda-\Lambda]1^T=:\frac{\alpha^2}{2}P(t,S).   \end{equation}

Let  \be\label{8} Q(t)=[I+(t\Gamma)^{-1}]^{-1} = I-(I+t\Gamma)^{-1}.\ee
For all sufficiently large $t$ the expression in (\ref{8}) renders itself to representation as an absolutely convergent geometric series:
\be
Q(t)= \sum_{v=0}^{\ff} (-1)^v\Big(\frac{\Ga^{-1}}{t}\Big)^v \label{9}
\ee
Hence for all $(i,j)$,
\be
q_{i,j}:= Q_{i,j} = \delta_{i,j} -\frac{(\Gamma^{-1})_{i,j}}{t} + \frac{(\Gamma^{-2})_{i,j}}{t^2} .... \label{10}
\ee
We proceed to find the series expansion for $\Phi_2(t,S)$ about $S=0$  in terms of  $Q(t)$. Henceforth the parameter $t$ will be suppressed in the expressions for  $Q(t)$ and  $\Phi_i(t,S)$ for $i=1,2$.

   To obtain (\ref{12}) consider (\ref{7}) and substituting    for $\Lambda$ using (\ref {4}) write:
   \be
  \Lambda \tilde{\Ga}\Lambda - \Lambda= t(I-S)[(t\Ga)^{-1}+I-S]^{-1}(I-S)-t(I-S) \label{15} 
  \ee
  Hence
\bea
P(S) &=& t1\cdot\{(I-S)[(t\Gamma)^{-1} + I-S]^{-1}(I-S)+ (S-I)\}\cdot1^T\nn\\
&=&t1\cdot\{(I-S)[Q^{-1}-S]^{-1}(I-S)+(S-I)\}\cdot1^T\nn
\\&=&t1\cdot\{(I-S)[I-QS]^{-1}Q(I-S) +(S-I)\}\cdot1^T \label{16}
 \eea

Since  $\det (QS) <1$ for all $t$ sufficiently large,
$$[I-QS]^{-1} = \sum_{m=0}^\infty { (QS)^m}.$$
Therefore 
\begin{eqnarray} \label{17}
&&(I-S)[I-QS]^{-1}Q(I-S)+ (S-I) \label{4.8}\\
&&\qquad=(I-S) \(\sum_{m=0}^\infty { (QS)^m}\) Q(I-S) + (S-I) \nonumber\\ 
&&\qquad=\(\sum_{m=0}^\infty { (QS)^m}\)Q -S\(\sum_{m=0}^\infty { (QS)^m}\)Q-\nonumber\\
&&\qquad\qquad\qquad  \(\sum_{m=1}^\infty { (QS)^m}\)+S\(\sum_{m=1}^\infty { (QS)^m}\)+ (S-I)\nonumber.
\end{eqnarray}
For $i=0,1$ we write 
$$
   S\(\sum_{m=i}^\infty { (QS)^m}\) = Q^{-1}\(\sum_{m=i+1}^\infty { (QS)^m}\), 
  $$ 
      to see that the last line of (\ref{17}) 
\[  =(Q-I)+(I-Q^{-1})\(\sum_{m=1}^\infty{(QS)^m}\)(Q-I)\nn  .
\]
This gives us (\ref{12}).

\qed

{\bf Proof of Theorem \ref{main1}}

 Let 
  \begin{equation} \label{1}
  D_i = \sum_{k=1}^n \Gamma^{-1}_{i,k} ; 
  \end{equation}
   i.e $D_i$ denotes the $i$-th row sum of $\Gamma^{-1}$.

By (\ref{8}):
\be I-Q^{-1} =I-(I+(t\Gamma)^{-1})=\frac{-\Gamma^{-1}}{t}, \label{19}\ee
 and \be Q-I=\frac{-\Gamma^{-1}}{t}(1 + O(1/t)).\label{20} \ee
Let $\overline{D}$ be a diagonal matrix such that $\overline{D}_{i,i}= D_i$.
Now substituting this in   (\ref{12}) we see that

\bea \label{21}
&&1\cdot[t (I- Q^{-1})\sum_{m=1}^\infty{(QS)^m}(Q-I)]\cdot1^T\\
&&\qquad  =\frac{1}{t}\sum_{m=1}^\infty1\cdot\Ga^{-1}(QS)^{m}\Ga^{-1 }(1+O(1/t))\cdot 1^T\nn\\
 &&\qquad=\frac{1}{t}\sum_{m=1}^\infty\sum_{i,j=1}^n \{\Ga^{-1}(QS)^{m}\Ga^{-1 }\}_{i,j}(1+O(1/t))\nn\\
 && \qquad=\frac{1}{t}\sum_{m=1}^\infty\sum_{i,j=1}^n \sum_{k,l=1}^n\Ga^{-1}_{i, k}(QS)^{m}_{k,l}\Ga^{-1 }_{l,j}(1+O(1/t))\nn\\
 &&\qquad=\frac{1}{t}\sum_{m=1}^\infty\sum_{k,l=1}^n(\sum_{i=1}^n \Ga^{-1}_{i, k})(QS)^{m}_{k,l}(\sum_{l=1}^n\Ga^{-1 }_{l,j})(1+O(1/t))\nn\\ 
&& \qquad=\frac{1}{t}\sum_{m=1}^\infty \sum_{k,l=1}^n(QS)^m_{k,l}D_kD_l(1+O(1/t))\nn
 \eea 
 
 \bea
&& \qquad=\frac{1}{t}\sum_{m=1}^\infty \sum_{k,l=1}^n \Big\{\overline{D}(QS)^m \overline{D} \Big\}_{k,l}(1+O(1/t))\nn\\
&& \qquad=\frac{1}{t}\sum_{k,l=1}^n \Big\{\overline{D}\(\sum_{m=1}^\infty (QS)^m \)\overline{D} \Big\}_{k,l}(1+O(1/t))\nn
 \eea
 
Now substitute this into (\ref{12}) to get (\ref{18}).  
\qed

 {\bf Proof of Lemma \ref{q pos}}  
   
 This is the direct consequence of the following
 
 \begin{lemma}\label{sign}
 Let $\Ga^{-1}$ be an $n$- dimensional $M$ matrix and assume that for some $i\ne j$, $\Gamma^{-1}_{i,j} =0$. Let $$k=\min\{l>1: \Ga^{-l}_{i,j}\neq 0\}$$
Then
\begin{itemize}
 \item[] \be k \leq n-1 \label{30}\ee
 \item[] and \be(\Gamma^{-k})_{i,j} = (-1)^k |(\Gamma^{-k})_{i,j}|.\label{31}\ee 
 
\end{itemize}
\end{lemma}

 \Proof
To show  (\ref{31}) make the following claim: if for all $1\leq u < k$, $\Gamma^{-u}_{i,j}=0$; then for any $h<k$,  any term of the form
$$\sum_{r_1}...\sum_{r_h} {\Gamma^{-1}_{i,r_1}\Gamma^{-1}_{r_1,r_2}....\Gamma^{-1}_{r_h,j}}$$ equals $0$.

For $k=2$; \begin{equation}(\Gamma^{-2})_{i,j} = \sum_{r=1}^n{\Gamma^{-1}_{i,r}\Gamma^{-1}_{r,j}}\geq 0\end{equation}
since all summands are the products of two negative factors, hence positive -- the only potential negative summands are of the form $\Gamma^{-1}_{i,i}\Gamma^{-1}_{i,j}$, which is zero here, hence, if $(\Gamma^{-2})_{i,j}=0$, then for all $1 \leq r\leq n$ ,
$\Gamma^{-1}_{i,r}\Gamma^{-1}_{j,r} =0$ -- so the above claim holds.

In general, for arbitrary $k$, suppose the claim holds for all $l<k$; $\Gamma^{-k}_{i,j}=0$. Then
 
\be
(\Gamma^{-k})_{i,j} = \sum_{r_1}\sum_{r_2}........\sum_{r_{k-1}}\Gamma^{-1}_{i,{r_1}}\Gamma^{-1}_{{r_2},{r_3}}........
\Gamma^{-1}_{{r_{k-1}}, j}.\label{32}
\ee
 Suppose a particular summand contains one or more factors of the form $\Gamma^{-1}_{i,i}$.\\ Rearranging the summation order we get:
\\
 \be
 \sum_{r_1}\sum_{r_2}...\sum_{r_u}\Gamma^{-1}_{r_1, r_1}\Gamma^{-1}_{r_2,r_2}...\Gamma^{-1}_{r_u,r_u}
\Big[\sum_{r_{u+1}}...\sum_{r_{k-1}}\Gamma^{-l}_{i,r_{u+1}}\Gamma^{-l}_{r_{u+1},r_{u+2}}...\Gamma^{-l}_{r_{k-1},j}\Big].\label{33}
\ee

However, each of the summands in square brackets is $0$ by assumption of our claim as it applies to $k-u <k$. So any potentially non-zero terms contain factors of the form $\Gamma^{-1}_{a,b}$ with $a\neq b$ only. But all such terms have the same sign as $(-1)^k$-- so that no cancellation is possible.  Hence $\Gamma^{-k}_{i,j}$ is either $0$ or has the sign of $(-1)^k$. This is (\ref{31}).
 \medskip

Let $E = \{E_{i,j}\}$ be the eigenvector matrix of $\Ga^{-1}$ and $U$ be its matrix of eigenvalues.

To obtain (\ref{30}) we claim that $(\Gamma^{-k})_{i,j}$ must be non-zero for some $k \leq n$; otherwise the eigenvectors of $\Gamma$ will have to satisfy a  homogeneous system of $n+1$ linear equations:
\begin{eqnarray} \sum_{k=1}^n{E_{i,k}E_{k,j}} \qquad =0\nonumber\\
\sum_{k=1}^n{E_{i,k}E_{k,j} u^{-1}_k} =0\nonumber\\
\sum_{k=1}^n{E_{i,k}E_{k,j} u^{-2}_k} =0\nonumber\\
.\nonumber\\
.\nonumber\\
\sum_{k=1}^n{E_{i,k}E_{k,j} u^{-(n-1)}_k }=0.\label{35}
\end{eqnarray}

 Let ${\mathbf{E}} = (E_{i,1}E_{1,j},E_{i,2}E_{2,j},....E_{i,n}E_{n,j})$ and write (\ref{35}) as: $$U \mathbf{E}^T={\mathbf{0}}$$ where $U$ is a square matrix with $U_{i,j}=u_j^{-{i-1}}$ for $1 \leq i,j \leq n$.

 We know that ${\mathbf{E}}$ is a non-zero vector -- otherwise $\Ga_{i,j}$ would be $0$. Hence $\det U$ must be $0$ -- that is it must have linearly dependent rows. Let  $$ \mathbf{u}= ( u_1^{-1}, u_2^{-1}, \ldots  u_n^{-1})$$
 and $$ \mathbf{u}^k= ( u_1^{-k}, u_2^{-k}, \ldots  u_n^{-k})$$
 Hence $ \mathbf{u}^{k-1}$ is the $k$-th row of $U$. The linear dependence of these vectors then implies that for some real numbers $a_k$; for $0 \leq k \leq n-1$,
 $$\sum _{k=0}^{n-1} {a_k \mathbf{u}^k}=\mathbf{0}.$$
 Which is, writing the vector form explicitly:
 For all $i$, $$ \sum _{k=0}^{n-1} {a_k u_i^k}=0.$$
 So each of eigenvectors of $\Ga^{-1}$ is the solution to the same polynomial  of degree $n-1$.

  By the Fundamental Theorem of Algebra that would mean that at least two of our eigenvectors are identical. WLOG, let those be $u_1$ and $u_2$. Let $\tilde{{\mathbf{E}}} = (E_{i,1}E_{1,j}+E_{i,2}E_{2,j},....E_{i,n}E_{n,j})$. Applying now the reasoning exactly identical to the above to the $n-1$ dimensional vector $\tilde{{\mathbf{E}}}$ and matrix $\tilde{U}$, such that $\tilde{U}_{i,j}=u_{j+1}^{-{i-1}}$ for $1\leq i,j \leq n-1$ conclude that $u_j=u_l$ for some $j \neq l$. Proceeding by induction conclude that all eigenvalues of $\Ga$ are equal. This contradicts the irreducibility of $\Ga$.
 
 Hence we know that (\ref{31}) is indeed non-zero for some  $k< n$. 
 
  \qed
  
  {\bf Proof of Lemma \ref{q pos} concluded.}
  
 \medskip 
  
Using this result to substitute into (\ref{33}) we conclude that whenever $\Gamma^{-l}_{i,j} =0$ $ \forall l<k$ and $\Gamma^{-k}_{i,j}\neq 0$,
 \begin{equation}
q_{i,j}= \frac{(-1)^k}{t^k}\Gamma^{-k}_{i,j} + O(t^{-(k+1)}) = \frac{1}{t^k}|\Gamma^{-k}_{i,j}| + O(t^{-(k+1)}) ,\label{34}
\end{equation}
 and so is positive for $t$ sufficiently large.
 
 \qed  
  
We will say that $\Gamma^{-1}$ has a zero of order $h$ at $(i,j)$ if $\forall r\leq h \\   \Gamma^{-r}_{i,j}=0$ and $\Gamma^{-(h+1)}_{i,j}\neq 0$.
\medskip
\begin{remark}
Note that the results of Lemmas \ref{div cond}, \ref{sign}, Theorem \ref{main1} and Corollary \ref{phi2} also hold for matrices that are not strictly positive definite.
\end{remark}
  
 \section{Upper bound on the critical point for Gaussian vectors.}

{\bf Proof of Theorem \ref{my thm} and Corollary \ref{0}.} 
 To obtain the upper bound on $\al_0$ from the Theorem \ref{main1} and Lemma \ref{q pos} we use  Lemma \ref{div cond}.
To apply the Lemma \ref{div cond} we extract the coefficient of a special type of term in the following

\begin{lemma}\label{coef}
 Let $A_{m_{i}, m_{j}}$ and  $B_{m_{i}, m_{j}}$ be the coefficients of the term $s_{i}^{m_{i}}s_{j}^{m_{j}}$ $(i \neq j)$ in $\Phi_{1}(S)$ and $\Phi_{2}(S)$ respectively. Then for all $ m_{i},m_{j}\ge 1$
\begin{equation}
 A_{m_i ,m_j } =\frac{1}{2}\Big(\frac{\Gamma^{-1}_{i,j}}{t}\Big)^2 (1+O(1/t))\label{22}
 \end{equation}
and
\begin{equation}
   B_{m_i ,m_j } = \frac{\alpha^2}{t^2} D_iD_j(-\Gamma^{-1}_{ i, j})(1+O(1/t))\label{23} 
   \end{equation}
   whenever $\Ga^{-1}_{i,j} \neq 0$.
 \end{lemma}

\Proof  Let $m=m_{i}+m_{j}$. It is clear  from (\ref{11}) that  the only contribution to  $ A_{m_i ,m_j } $ comes from 
\begin{equation}
\frac{ \{ (QS)^{m}\}_{i,i}+\{ (QS)^{m}\}_{j,j}}{  m},\label{24}
   \end{equation}

and, furthermore,  we must have 
\begin{equation}
 \frac{  \{(QS)^{m}\}_{i,i}}{ s_{i}^{m_{i}}s_{j}^{m_{j}}}=\prod_{l=1}^{m}q_{p_{0},p_{1}}q_{p_{1},p_{2}}\cdots q_{p_{l-1},p_{l}}\cdots q_{p_{m-1},p_{m}} \label{25}
   \end{equation}
with $p_{0}=p_{m}=i$; and all the other $p_{l}$, $1\le l\le m-1 $, must be either $i$ or $j$. Since
\begin{equation}
   q_{i,i}=1+O(1/t),\quad q_{j,j}=1+O(1/t)\quad\mbox{and}\quad q_{i,j}=q_{j,i}=-\frac{\Ga^{-1}_{i,j}}{t}+O(1/t^{2}),\label{26}
   \end{equation}
the terms on the right-hand side of (\ref{25}) that are not $O(1/t^{3})$ are those terms in which $q_{i,j}$ and $q_{j,i}$ each occur only once. This can happen in the following  $m_{i}$ ways
\begin{equation}
   q_{i,i}^{r}q_{i,j}q_{j,j}^{s}q_{j,i}  q_{i,i}^{u},\qquad r=0,\ldots,m_{i}
   \end{equation}
(and, obviously, $s=m_{j}-1$ and $u=m_{i}-r $). Clearly
\begin{equation}
      q_{i,i}^{r}q_{i,j}q_{j,j}^{s}q_{j,i}  q_{i,i}^{u}=\frac{(\Ga^{-1}_{i,j})^{2}}{t^{2}}+O(1/t^{3}).
   \end{equation}

Repeating this argument with $i$ and $j$ interchanged we see that there are $m_{j}$ ways that 
\begin{equation}
    \frac{  \{(QS)^{m}\}_{j,j}}{ s_{i}^{m_{i}}s_{j}^{m_{j}}}=\frac{(\Ga^{-1}_{i,j})^{2}}{t^{2}}+O(1/t^{3}). \label{27}
   \end{equation}
 Using (\ref{25})--(\ref{27}) we get (\ref{22})
\medskip
It is clear  from (\ref{18}) that  the only contribution to the leading term of $ B_{m_i ,m_j } $ comes from 
\begin{equation} \label {27a}
\frac{\al^{2}}{2t }\(\sum_{u=1}^n  D_{u}\{ q_{u,i}s_i(QS)^{m-1}\}_{i,j}D_{j}+\sum_{u=1}^n  D_{u}\{ q_{u,j}s_j\{ (QS)^{m-1}\}_{j,i}D_{i}\);
   \end{equation}
   
   By (\ref{26})  we see that (\ref{27a}) equals to
   \be \label{27b}
 \frac{\al^{2}}{2t }\(D_i\{ (QS)^{m}\}_{i,j}D_{j}+  D_{j}\{ \{ (QS)^{m}\}_{j,i}D_{i}\)(1+O(1/t));
   \end{equation}  
where the contribution to  the first summand of   $ B_{m_i ,m_j } $ will be from  a coefficient arising from $m_i$ of $s_i$ factors followed by $m_j$ of $s_j$ factors; and to the the second --  from $m_j$ of $s_j$ factors followed by $m_i$ of $s_i$ factors.
Furthermore, by (\ref{26}) the part of  $ B_{m_i ,m_j } $ arising from $\{(QS)^{m}\}_{i,j}$ equals to
\begin{equation}
    q_{i,i}^{m_{i}-1}q_{i,j}q_{j,j}^{m_{j}}+O(1/t^{2})
   \end{equation}
and similarly with $i$ and $j$ interchanged. Since
\begin{equation}
   q_{i,i}^{m_{i}-1}q_{i,j}q_{j,j}^{m_{j}}=\frac{-\Ga^{-1}_{i,j}}{t}+O(1/t^{2})
   \end{equation}
   we get (\ref{23}). \qed

\noindent{\bf Proof of Theorem \ref{my thm} continued.}

To begin, assume that $\Ga^{-1}_{i,j} \neq 0$ for all $(i,j)$.

Let 
 $$C_{m_i, m_j}=A_{m_i, m_j}+B_{m_i, m_j}$$
From (22) and (23):

\be
C_{i,j}(\al):=C_{m_i, m_j}= \frac{-\Gamma^{-1}_{i,j}}{t^2}\(\frac{-\Gamma^{-1}_{i,j}}{2} +\alpha^2 D_iD_j\)(1+O(1/t))\label{28}
\ee
Note that since the length of the term is fixed, for any $N$ such that $m_i, m_j <N$ (\ref{28}) holds for all $t$ sufficiently large (i.e. $t>>N$) thus making $C_{i,j}(\al)$ independent of   $m_i, m_j$.

For $C_{i,j}(\alpha)$ to remain non-negative when $D_iD_j <0$ as $t\to\infty$ we need to have 
 
\begin{equation}
\frac{|\Gamma^{-1}_{i,j}|}{2}+\alpha^2D_iD_j >0.\label{29}
\end{equation}
Therefore, using Lemma \ref{div cond}, we get (\ref{3}), which is Theorem \ref{my thm} for  inverse covariance matrices with non-zero entries.
\medskip 
 
 \noindent{\bf Removing restriction  $\Ga^{-1} \neq 0$ }

Now returning to (\ref{22}) and (\ref{23}) we see that if $\Gamma^{-1}$ has a zero of order $k$ at $(i,j)$ for some $0<k< n-1$ then using Theorem \ref{main1} again we obtain:
\begin{equation}
A_{i,j}= \frac{1}{2t^{2(k+1)}}((\Gamma^{-k})_{i,j})^2 + O(t^{-(2k+3)}) \label{37}
\end{equation}
And:\begin{equation}
B_{i,j}= \alpha^2 D_iD_j \frac{1}{t^{k+2}}|(\Gamma^{-k})_{i,j}| +O(t^{-(k+3)})\end{equation} \label{38}
So that
\begin {equation}\label{39}
C_{i,j}(\al) =\frac{1}{2t^{2(k+1)}}((\Gamma^{-k})_{i,j})^2 +\alpha^2 D_iD_j \frac{1}{t^{k+2}}|(\Gamma^{-k})_{i,j}| (1+O(1/t))\\
\end{equation}
\medskip

 will be positive iff $D_{i}D_{j} >0$, i.e. whenever $i$-th and $j$-th rows have the same sign. This concludes the proof of Theorem \ref{my thm} and Corollary \ref{0}.
\qed

\section{ Elementary proof of Theorem \ref{var mean}}

We are now in the position to give the "`elementary"' - i.e. not involving the Markov process theory - proof of the Theorem \ref{var mean}.This proof simplifies  the proofs  given in \cite { inf div}, Theorem 1.3; see also \cite{book}, Theorem 13.3.1. We will start by showing the result for the special case $c=\mathbf{1}$, see \cite{inf div}, Theorem 1.2:

\begin{theorem}\label{unif mean}
  Let $G$ be a mean zero Gaussian vector with strictly positive definite covariance matrix $\Gamma$. The following are equivalent:
 \begin{enumerate}\item$(G+\al)$ has infinitely divisible squares for all $\al \in {\mathbf{R}}$. 
 
\item $\Gamma^{-1}$ is an $M$-matrix with non-negative row sums. 
 \end{enumerate} 
  \end{theorem}

{\bf Proof}

$(1 \Rightarrow 2)$: Let $(G+\al)$ have infinitely divisible squares for all $\al \in {\mathbf{R}}$. By Theorem 13.2.1, \cite{book} $\Ga^{-1}$ is an $M$-matrix. By Corollary \label{0, infty} $\DD = \emptyset$ in this case. Since $\Ga$, and so also $\Ga^{-1}$ is positive definite this implies that all row sums of $\Ga^{-1}$ are non-negative. 

$(2 \Rightarrow 1)$: Let $\Gamma^{-1}$ be an $M$-matrix with non-negative row sums. Its Laplace Transform is given by Theorem \ref{main1}. The positivity of coefficients in $\Phi_1$ follows from $\Gamma^{-1}$ being an $M$-matrix; and that of coefficients in $\Phi_2$ from Lemma \ref{q pos} and row sums positivity.

{\bf Proof of Theorem \ref{var mean}}

We first observe that if $ (1), (2) or(3)$ of Theorem \ref{var mean} hold then
\begin{itemize}
\item{} All components $c_i$ of vector $c$ must have the same sign.
\item{} If $c_i = 0$ for some $i$, then $c=\bf {0}$.
\end{itemize}

\begin{itemize}

\item[] If $(1)$ holds then no component of $c$ is $0$ or $c= \bf{0}$.
Suppose $c_i=0$,  $c_j \neq 0$ or $c_i, c_j$ have different signs for some pair $(i,j)$. Since $(\Ga + \al c)^2$ is infinitely divisible for all $\al$, so is the sub vector  $( \Ga_i+\al c_i, \Ga_j+ \al c_j)$. But for this sub vector $\Ga^{-1} C$has row sums of different sign. Hence it is not infinitely divisible for all $\al$ by Remark\ref{ub}. Hence $c_i \neq 0$ for all $i$ and all components of $c$ are of the same sign.

\item[] If $(2)$ holds then it also holds for any sub vector $(\Ga_i + c_i \eta, \Ga_j + c_j \eta, \eta)$. However by the same reasoning as above either $\Ga^{-1}_{i,3}$ or $\Ga^{-1}_{j,3}$ is negative by the Lemma 2.1, \cite {inf div}.
\item[] If $(3)$ holds then recalling that $\Ga^{-1}$ is strictly positive definite $M$-matrix conclude that $c_i$'s are non-zero and of the same sign.
\end{itemize} 

The proof of Theorem 1.3 in \cite {inf div} now shows that the result of \ref{var mean} follows immediately by replacing $G$ by $\frac{G}{c}$.

 \section{Some other applications.}  
   
Another open question pertaining to infinite divisibility of the Gaussian squares is: Suppose $G^2$ is not infinitely divisible. What could be said about infinite divisibility of $(G+{\mathbf c})^2$? By the technique employed in Theorem \ref{var mean} it is enough to consider this question for $(G+\al {\mathbf 1})^2$ for $\al \in  {\mathbf R}$. Furthermore, by conjugating $G$ by an appropriate signature matrix one can consider vector $G$ with covariance $\Ga$ having all non-negative entries. Since in two dimensions the inverse of covariance matrix is always an $M$-matrix, one should be concerned with vectors in three or more dimensions.

Hence the question stands as follows: consider the gaussian vector $G$ of dimension 3 or higher and assume that $G^2$ is not infinitely divisible. Let $\Ga$ be the covariance matrix of $G$ such that $\Ga_{i,j} \geq 0$ for all$(i,j)$. When is $(G+\al)^2$ infinitely divisible for all $\al$? For no values of $\al$? Are there some vectors possessing the critical point - i.e. some $\al_0 >0$ such that  $(G+\al)^2$ is not infinitely divisible for $\al \leq \al_0$ and infinitely divisible for $\al > \al_0$?

Theorem \ref{main1} allows us to shed some light on the above question for some special cases.
\bl \label{non id}

Let $G$ be a Gaussian vector with positive, strictly positive definite covariance $\Ga$. 
If there exist a triple $i,j,k$ such that 
\begin{eqnarray}
\label{0.1}
\Ga^{-1}_{k,j}, \Ga^{-1}_{i,k} <0 \quad \mbox{and} \quad \Ga^{-1}_{j,i} >0,\\
\mbox{and}\quad  D_i, D_j, D_k >0
\end{eqnarray}
 then

$(G+\al)^2$ is not infinitely divisible for any $\al \in \mathbf{R}$.
\el
\Proof

Consider the coefficient of the term $s_is_j$:
\begin{eqnarray}
&&\frac{1}{2}\big(q_{i,j}^2+\alpha_i^2\frac{2D_iD_j}{t}q_{i,j}\big)
(1+O(1/t))\nonumber\\
=&& \frac{1}{2t^2}\big((\Gamma_{i,j}^{-1})^2-\alpha|D_iD_j|\Gamma_{i,j}^{-1}\big) (1+O(1/t))
\end{eqnarray}
This is negative when
\begin{equation}
\alpha^2 > \frac {\Gamma_{i,j}^{-1}}{2D_iD_j}
\end{equation}

Now consider the term $s_is_js_k$. Its coefficient is:
\begin{eqnarray}
&& \big\{q_{i,j}q_{j,k}q_{k,i} + \frac{\alpha^2}{t}[|D_iD_j|q_{j,k}q_{k,i}+|D_jD_k| q_{i,j}q_{k,i}+|D_kD_i|q_{i,j}q_{j,k}]\big\}\nn\\
&&\qquad\qquad\qquad\qquad\qquad\qquad\qquad\qquad\qquad\qquad\qquad \cd (1+O(1/t))\nonumber\\
=&&\frac {1}{t^3} \big\{ -\Gamma^{-1}_{i,j}\Gamma^{-1}_{j,k}\Gamma^{-1}_{k,i}+ \nn\\&& \alpha^2[D_iD_j|\Gamma^{-1}_{j,k}\Gamma^{-1}_{k,i}| -D_jD_k |\Gamma^{-1}_{i,j}\Gamma^{-1}_{k,i}| -D_kD_i|\Gamma^{-1}_{i,j}\Gamma^{-1}_{j,k}|]\big\}(1+O(1/t))\nonumber\\
\end{eqnarray}
Note that only the first term in square brackets is positive, while the other three terms are negative. Hence for this expression to be negative it is enough to have:
\begin{equation}
\alpha^2 < \frac{\Gamma^{-1}_{i,j}}{D_iD_j}.
\end{equation}

This shows that the series expansion of $\Phi$ will have at least one negative (non-constant term) coefficient for each $\al \in {\mathbf R}$.

\qed
\begin{corollary} \label{dim 3,4,5}
Let $n=3, 4 \mbox{ or } 5$. Assume that  $\Ga >0$,  is strictly positive definite and $\Ga^{-1}$ has positive row sums. If $\Ga^{-1}$ has at least one positive off-diagonal entry, then
$(G+\al)^2$ is not infinitely divisible for all $\al \in \mathbf{R}$.
\end{corollary}

In general case of $\Gamma$, there exist a signature matrix $N$ such that $\tilde{\Gamma} = N\Gamma N >0$ and replace $\Ga$ by $\tilde {\Ga}$ throughout.

\Proof
Note that to obtain this result it is enough to show that in dimension 3, 4 and 5 the inverse covariance matrix with positive off-diagonal elements has a triple of indices satisfying (\ref{0.1}). Also note that $\Ga ^{-1}$ can not have a row with only positive entries.

This is obvious for $n=3$, since in this case one off diagonal element is positive, while two others are negative.

To see the result for $n=4 \quad \mbox{and}\quad 5$ assume that $\Ga^{-1}$ does not contain the required combination of entries. Note that in this case  it is possible by means of row-column permutations to bring the matrix $\Ga^{-1}$ into the following form:
$$ \mbox{if  }\qquad j>i\qquad \mbox{and} \qquad\Ga^{-1}_{i,j} >0, \mbox { then }\Ga^{-1}_{i,k} >0, \mbox{ for all } k>j$$
When written in this form, it is easy to see that we must have $\Ga^{-1} = A+B$, where
$A$ is the direct sum of $M$-matrices;  $B \geq 0$ whenever $A=0$, and $B=0$ whenever $A \neq 0$.
If $n=4 \mbox{ or } 5$ this is only possible if $A$ is a direct sum of two such matrices, each of the size $2 \times 2$ for $n=4$, or of size $2\times 2$ and $3 \times 3$ for $n=5$. However, if this is the case, then we can find a signature matrix $N$ such that $N\Ga^{-1}N$ is a weakly diagonally dominant matrix with negative off diagonal and positive diagonal elements. But the by \cite{book}, Remark 13.1.3, $(N\Ga^{-1}N)^{-1}= N\Ga N$ has  positive entries. Contradiction. Hence the  triple of entries satisfying \ref{0.1} indeed exists for $n=4 \quad \mbox {and}\quad 5$.
\qed

 \section{General form of the coefficient for the term of the form $s_i^{m _i}s_j^{m_j}$.}
 Theorem \ref{main1} can be used to derive the two index coefficient of $\Phi$ explicitly for any dimension $n$. When applied to the two-dimensional case this result gives a shortcut for the derivation of expression for $R_{i,j,p}$ used by Marcus and Rosen in \cite{critical}.
 
In the statement of Lemma \ref{coef} we define  $A_{m_{i}, m_{j}}$ and  $B_{m_{i}, m_{j}}$ to be the coefficients of the term $s_{i}^{m_{i}}s_{j}^{m_{j}}$;  $(i \neq j)$ in $\Phi_{1}(S)$ and $\Phi_{2}(S)$ respectively. In (\ref{22})  we estimate $A_{m_{i},m_{j}}$, for fixed $m_{i}$, $m_{j}$ as $t\to \ff$. However, to apply Lemma \ref{div cond} to show coefficient positivity we must also consider $A_{m_{i},m_{j}}$ when $m_{i}$ and $m_{j}$ are larger than $t$.
We do this in the next lemma.

\begin{lemma} For $m_i, m_j \geq 1$ 
\begin{eqnarray}
A_{m_i, m_j}&= &\sum_{p=0}^{min(m_i, m_j)-1}\binom{m_i-1}{p}\binom{m_j-1}{p}\label{2.1}\\ 
 &&  \qquad  \qquad \cd \frac{1}{p+1}(q_{i,i})^{m_i-p-1}(q_{j,j})^{m_j-p-1}(q_{i,j})^{2(p+1)}\nn\\
& =&\sum_{p=0}^{min(m_i, m_j)-1}V(m_i,m_j,p)\frac{q_{i,j}^{2}}{p+1} 
 (1+O(1/t))\nn,
\end{eqnarray}
where
\be
V(m_i,m_j,p) = (q_{i,i}^{m_i-(p+1)})(q_{j,j}^{m_j-(p+1)})(q_{i,j}^{2p})\binom{m_i-1}{p}\binom{m_j-1}{p}.
\ee
 \end{lemma}
 
 \Proof 
By (\ref{11}), to find $A_{ m_i, m_j}$ we   sum over all  terms in the trace of $(QS)^{m_i+m_j}$ that involve $s_i$, $m_i$ times and $s_j$, $m_j$ times. Since the summands of the coefficient come from the trace of a matrix each of them will have the form $q_{i_1,i_2}q_{i_2,i_3}....q_{i_{m-1},i_m}$ with $i_1=i_m=i$  or $i_1=i_m=j$.
 It is  convenient to consider the circular arrangement of the  factors $s_{(\cd)}$, since once such arrangement is fixed each of  $m$ of its rotations   result in the same coefficient.

Consider arrangement of the $s_i$ and $s_{j}$ into blocks. By   blocks we mean an unbroken string of $s_i$ or $s_j$ factors of   length one or greater. For example, $s_{j}^{3}s_i^4 s_{j}$ contains an $s_{i}$ block of  length 4  but $s_is_j^2s_i^3$ does not.

To account for all possible arrangements of indices of $s_i^{m _i}s_j^{m_j}$,  consider any fixed circular arrangement. Suppose this arrangement has $p+1$ separate blocks of $s_{i}$,  it then must also have $p+1$ blocks of $s_{j}$. We will have $0\leq p < \min\{m_i, m_j\}$. For each $p$  there exist 
\be
\binom {m_i-1}{p}\binom {m_j-1}{p}
\ee
ways to make this separation into blocks. Once the blocks are defined there exist exactly one circular arrangement of blocks( if we agree which index to always place first).

Each arrangement can be rotated one position clockwise exactly  $m=m_i+m_j$ times. A given arrangement will repeat itself after $k$ one-positional rotations, where $m=kv_1$ for some integer $v_1$. Obviously also $p+1= kv_2$  for some integer $v_2$. Hence this arrangement can be rotated $v_1$ times giving rise to a new arrangement. This cancels the $v_1$ factor in the denominator of (\ref{11}), leaving $k$. The coefficients in $\Phi_1$ of all such rotations are identical and equal to $$\frac{1}{m}\{q_{i,j}^{p+1}q_{i,i}^{m_i-p-1}q_{j,j}^{m_j-p-1}\}.$$
To see this note that after any such rotation the resulting string of factors will either begin and end with the same index or will begin with $s_i$, end with $s_j$ or vice versa. Using the symmetry of $Q$, see that:

If it begins and ends with, say, $s_i$, the coefficient in $\mbox{trace}(QS)^m$ is:
$$q_{i,i}(q_{i,i})^{m_i-p-2}(q_{j,j})^{m_j-p-1}(q_{i,j})^{p+1},$$
and is it begins with, say, $s_i$ and ends with$s_j$, the coefficient in $\mbox{trace}(QS)^m$ is:
$$q_{j,i}(q_{i,i})^{m_i-p-1}(q_{j,j})^{m_j-p-1}(q_{i,j})^{p},$$
and the same with $i, j$ interchanged. Using the fact that $q_{i,j}=q_{j,i}$ again, see that this equals to  the expression in figure parentheses above,although the order of $q$ factors varies with the rotation; for example starting with an arrangement giving rise to a coefficient 
$$\frac{1}{m} q_{i,i}^{u_1}q_{i,j}q_{j,j}^{v_1}.....q_{i,j}q_{j,j}^{v_{p+1}}q_{j,i}$$
rotating it one position results in a coefficient
$$\frac{1}{m} q_{j,i}q_{i,i}^{u_1}q_{i,j}q_{j,j}^{v_1}.....q_{i,j}q_{j,j}^{v_{p+1}}$$
and so on. 
Also note that each fixed arrangement with a fixed starting point gives rise to two coefficient summands, each corresponding to the direction - clockwise or counterclockwise - of the order of $q$-factors. For example, the last arrangement "`read"' in the opposite direction will give the coefficient
$$\frac{1}{m} q_{j,j}^{v_{p+1}}q_{j,i}..... q_{j,j}^{v_1}q_{j,i}q_{i,i}^{u_1} q_{i,j}.$$
Therefore the coefficient of the whole sum is multiplied by $2$.

 Finally we note that each   of the $v_1$ rotated block arrangements give rise to $v_2$ repetitions, according to which block is taken as the starting point   of the circular arrangement.  Therefore  we divide by $v_2$. Thus the resulting denominator is: $v_2 k=p+1$. \qed
 
 \begin{remark}
 We pause in the development of the two-index coefficients to remark that the very argument just employed to find the denominator of $A_{m_{i},m_{j}}$ part of the coefficient is responsible for the non-existence of the nice "`closed"' form for the coefficient in tree or more indices.
 Consider the coefficient of  $s_i^{m _i}s_j^{m_j}s_k^{m_k}$ in the power series expansion of $\Phi_{1}(S)$, and assume that $m_i, m_j, m_k >0$. Let our factors be broken into $p_i, p_j, p_k$ blocks respectively and assume an circular arrangement is fixed such that no two blocks of the same index are adjacent.  We can rotate this arrangement $u$ times to get a new arrangement, and obviously $m=m_i+m_j+m_k = uv_1$ for some integer $v_1$. Clearly, too, $p_i= uv_i, p_j=uv_j, p_k=uv_k$. 
 Now divide by $m$ and multiply by $u$ to get $v_1$ in denominator as before. Now to account for different starting blocks we would have to divide by $v_i, v_j$ or $v_k$, depending on which block sequence gives $u$ new arrangements under rotation. But to know this we need additional information on the blocks breakdown! This is what will prevent us from getting the sought coefficient based on the number of blocks alone.
 \end{remark}

Recall that $B_{m_{i},m_{j}}$ is the coefficient of  $s_i^{m _i}s_j^{m_j}$ in the power series expansion of $\Phi_{2}(S)$. Let
$$B_{m_{i},m_{j}} = \sum_{p=0}^{min(m_i, m_j)-1}B_{m_i,m_j,p}$$

\begin{lemma}  For $m_i, m_j \geq 1$
\begin{eqnarray} \label{2.6}
 B_{m_i,m_j}
& = &\sum_{p=0}^{min(m_i, m_j)-1}\frac{\al^2}{2t} V(m_i,m_j,p)\big\{2D_iD_j \nn\\
 &&\qquad +\big [\frac{p}{m_i-p}D_j^2+\frac{p}{m_j-p}D_i^2\big]\big\}(1+O(1/t)).
\end{eqnarray}
\end{lemma}

 \Proof
By  (\ref{18}), (\ref{7}), (\ref{19}) and (\ref{20}) we obtain the coefficient in $P(S)$ arising from a particular arrangement, say $s_{k_1}s_{k_2}s_{k_3}...s_{k_{m_i+m_j}}$, where $k_u$ is either $i$ or $j$ equals to:
\begin{eqnarray} \label{2.20}
&&\frac{1}{t}\sum_{r=1}^n D_{r}q_{(r,k_1)}q_{(k_1,k_2)}q_{(k_2, k_3)}\cdots\nn\\
&&\qquad\cdots q_{(k_{m_i+m_j-1}, k_{m_i+m_j})}q_{(k_{m_i+m_j}, k_{m_i+m_j})}D_{k_{m_i+m_j}}(1+O(1/t))\nn 
\end{eqnarray}

 Using  \ref{26} see that  (\ref{2.20}) equals to:
\begin{eqnarray}
&&\frac{1}{t}D_{k_1}q_{(k_1, k_2)}q_{(k_2, k_3)}\cdots\nn\\
&&\qquad\cdots q_{(k_{m_i+m_j-1}, k_{m_i+m_j})}q_{(k_{m_i+m_j}, k_{m_i+m_j})}D_{k_{m_i+m_j}}(1+O(1/t))\nn 
\end{eqnarray}

To find  $B_{m_i,m_j,p}$ , assume that one  group of factors - either $s_i$ or $s_j$ - is broken into $p+1$ blocks for some $p$ . Since we are dealing with a term in two indices only, one of the two cases is possible:
\begin{enumerate}
\item{} Both groups of factors are broken into $p+1$ blocks. Then $0 \leq p < \min\{m_i, m_j\}$. There will be $\binom{m_i-1}{p}\binom{m_j-1}{p}$ such arrangements. That means either starting with $s_i$, ending with $s_j$ block, or vice versa. Both will result in exactly $2p+1$ of $q_{i,j}$ factors  and  have the  coefficient in $\Phi_2$:
\begin{eqnarray}\label{2.2}
&&\sum_{p=0}^{min(m_i, m_j)-1}\binom{m_i-1}{p}\binom{m_j-1}{p}\frac{\alpha^2}{2t}\Big[ D_i(q_{i,i}^{m_i}q_{j,j}^{m_j}) \frac{q_{i,j}^{(2p+1)}}{ (q_{i,i}q_{j,j})^{(p+1)}}q_{j,j}D_j\nn\\
&&\qquad+D_j(q_{i,i}^{m_i}q_{j,j}^{m_j})\frac{q_{i,j}^{(2p+1)}}{ (q_{i,i}q_{j,j})^{(p+1)}}q_{i,i}D_i\Big](1+O(1/t))\nn\\
&&=\sum_{p=0}^{min(m_i, m_j)-1}\frac{\alpha^2}{2t} D_iD_jV(m_i,m_j,p)(\frac{q_{i,j}}{q_{i,i}}+\frac{q_{i,j}}{q_{j,j}})(1+O(1/t))
\end{eqnarray}

\item{}Now restrict $p$ to $1 \leq p < \min\{m_i, m_j\}$. The group of $s_i$ factors is broken into $p+1$ blocks and that of $s_j$ factors into $p$ blocks, or vice versa. (If $m_i <m_j$ then there is also a case when there are $m_i$ \quad of $s_i$ factors and  $m_i+1$ \quad of $s_j$ factors. The resulting coefficient in $\Phi_2$ is of the form $B_{m_i,m_j,m_i}(1/t)$.) There will be $\binom{m_i-1}{p}\binom{m_j-1}{p-1}$ arrangements of the first kind and $\binom{m_i-1}{p-1}\binom{m_j-1}{p}$ of the second. This will result in exactly $2p$ of  $q_{i,j}$ factors, moreover the $D$  factors will have the same index in each of the arrangements - same index as the factors broken into $p+1$ blocks. The corresponding coefficient in $\Phi_2$ is, then:

\bea \label{2.3}
&&\frac{\alpha^2}{2t}\sum_{p=1}^{min(m_i, m_j)-1}(q_{i,i}^{m_i}q_{j,j}^{m_j})q_{i,j}^{2p}\nn\\
&&\qquad\qquad \cdot\Big[\binom{m_i-1}{p}\binom{m_j-1}{p-1} \frac{D_i^2 q_{i,i}}{ (q_{i,i})^{(p+1)}(q_{j,j})^{p}}\nn\\
&&\qquad+\binom{m_i-1}{p-1}\binom{m_j-1}{p}\frac{D_j^2q_{j,j}}{ (q_{i,i})^p (q_{j,j})^{(p+1)}}\Big](1+O(1/t))\nn\\
\eea
Using now the identity:
$$\binom{m-1}{p-1}= \frac{p}{m-p} \binom{m-1}{p}$$  see that (\ref{2.3}) equals to:
\bea \label{2.4}
\sum_{p=1}^{min(m_i, m_j)-1}\!\!\!\frac{\alpha^2}{2t} V(m_i,m_j,p)\big[\frac{p}{m_i-p}D_j^2+\frac{p}{m_j-p}D_i^2\big](1+O(1/t))\nn\\
\eea
since the expression in square parentheses is $0$ when $p=0$.
\end{enumerate}

Now denote as follows:
\begin{eqnarray}
A_{m_i, m_j,p}= V(m_i,m_j,p)\frac{q_{i,j}^2}{p+1} (1+O(1/t)) \label{2.5}
\end{eqnarray}
 
 And
\begin{eqnarray} \label{2.6}
 &&B_{m_i,m_j,p}
 = \frac{\alpha^2}{2t} V(m_i,m_j,p)\big\{\big[\frac{p}{m_i-p}D_j^2+\frac{p}{m_j-p}D_i^2\big]\nn\\
 &&\qquad\qquad\qquad\qquad\qquad\qquad+2D_iD_j\big\}(1+O(1/t))
\end{eqnarray}
 
 Let
 \bea
 C_{m_i,m_j, p}=A_{m_i, m_j,p}+ B_{m_i,m_j,p}
 \eea
 
 Hence,
 \be
C_{m_i,m_j}= \sum_{p=0}^{\min(m_i,m_j)-1} C_{m_i,m_j, p}
\ee
And
 
\bea
C_{m_i,m_j, p}=&&\quad  V(m_i,m_j,p)\Big\{ \frac{(q_{i,j})^2}{p+1}  \\
 &&\quad+  \frac{\alpha^2}{2t}\Big[2D_iD_jq_{i,j}
  + \frac{p}{m_j-p} D_i^2 +\frac{p}{m_i-p} D_j^2  \Big]\Big\}(1+O(1/t))\nn \label{2.7}
\eea
\bigskip

Denote the expression in figure parentheses by $R_{m_i,m_j, p}$. We are interested in the sign of $C_{m_i,m_j}$, which will coincide with the sign of $C_{m_i,m_j, p}$, if they are identical for all $p$. Since the sign of  $C_{m_i,m_j, p}$ is the same as the sign of $R_{m_i,m_j, p}$ we will concentrate on the later.
Assume that $m_im_j$ is bounded by linear multiple of $t^2$, i.e.
$$\sqrt{m_im_j} \leq Nt  \quad \mbox{for some fixed}\quad  N $$

\bea \label{2.8}
R_{m_i,m_j, p}&=&\quad \Big\{ \frac{1}{p+1} \frac{|\Ga^{-1}_{i,j}|^2}{t^2} +\frac{\alpha^2}{2} \Big[-\frac{2|D_iD_j||\Ga^{-1}_{i,j}|}{t^2}\nn\\
& &\qquad\qquad\qquad+ \big(\frac{D_i^2}{t}\frac{p}{m_j-p}+\frac{D_j^2}{t}\frac{p}{m_i-p} \big)\Big]\Big\}(1+O(1/t))\nn \\
 &\geq&\quad \Big\{ \frac{1}{p+1} \frac{|\Ga^{-1}_{i,j}|^2}{t^2} +\frac{\alpha^2}{2} \Big[-\frac{2|D_iD_j||\Ga^{-1}_{i,j}|}{t^2}\nn\\
& &\qquad\qquad\qquad+\frac{2|D_iD_j|}{t}\frac{p}{\sqrt{m_im_j}} \Big]\Big\}(1+O(1/t))\nn \\
&\geq &\quad \frac{1}{t^2}\Big[ \frac{|\Ga^{-1}_{i,j}|^2}{p+1}  
+\frac{\alpha^2}{2} |2D_iD_j|(-|\Ga^{-1}_{i,j}|+ \frac{p}{N} )\Big](1+O(1/t))\nn \\
\eea

The representation in ( \ref{2.8}) easily leads to the proof of Lemmas 5.2 and 5.3 of \cite{critical} by providing a combinatorial alternative to algebraic derivation of the expression for $R_{j,k,p}$ employed in that paper.

\section{Some numerical examples.}
 The Gaussian vector is called an assosiated with a Borel process when the covariance matrix is also a $0$-potential density of Borel right process.  The following (necessary but not sufficient) property of an associated Gaussian vector is the direct consequence of the property of $0$-potential density of Borel right process:
\be
\Ga_{i,j} \leq \Ga_{i,i} \land \Ga_{j,j}.\label{3.1}
\ee
contrasting with the regular property of the strictly positive definite covariance matrix:
\be
\Ga_{i,j}^2 \leq \Ga_{i,i}\Ga_{j,j}.\label{3.2}
\ee

Although the above regularity property is too weak to provide much information about $G$ in general; in case $n=3$ we have the following 
\begin{corollary}\label{3 dim}
 Let $G=(G_1, G_2, G_3)$ be a Gaussian vector satisfying (\ref{3.1}). Then the upper bound for the critical point of $G$ is non-zero.
 \end{corollary}
 \Proof
 This follows immediately from Theorem \ref{my thm} by writing the matrix $\Ga$ in terms of $\Ga^{-1}$ and above mentioned regularity condition:
 
 Suppose $\Gamma^{-1}$ has a zero entry $\Gamma^{-1}_{i,j}$ while $\Gamma$ satisfies (\ref{3.1}). Let $|\Ga^{-1}_{i,j}|= g_{i,j}$ and let $g_{1,3}=0$. Then $|\Ga_{1,2}|=g_{1,2}g_{3,3}$; $|\Ga_{2,3}|=g_{1,1}g_{2,3}$; $|\Ga_{2,2}|=g_{1,1}g_{3,3}$. This implies that $g_{1,2}<g_{1,1}$ and $g_{2,3}<g_{3,3}$. Hence $D_1$ and $D_3$ are both positive, which makes the first case of Corollary \ref{0, infty} impossible for such matrices.

The $G$-vectors with inverse covariance matrices $\Gamma^{-1}$ with zero entries and corresponding row sums of different signs indeed exist.\\
Below is an example for $n =3$ with one zero entry:

Let
 \begin{equation*} \Gamma = \frac{1}{7}\left(\begin{array}{ccc}
15&4&2\\
4&2&1\\
2&1&4\end{array}
\right);
\mbox{ then}\quad  
   \Gamma^{-1} = \left(\begin{array}{ccc}
1&-2&0\\

-2&8&-1\\

0&-1&2\end{array}
\right )\\
\end{equation*}

Where $\Gamma^{-1}$ is an $M$-matrix, so $G^2$ is infinitely divisible and 
       $\Gamma^{-2}_{1,3} = 2$ and  $\Gamma^{-1}_{1,3} = 0$, $D_1 = -1 <0$, $D_3 = 1 >0$.
        
Therefore for a vector with this covariance matrix $\tilde{B}_{i,j}=  -2\frac{1}{t^3}$   while  $\tilde{A}_{i,j}= \frac{1}{t^4}$ 
 so that $\tilde{C}_{i,j}=\frac{1}{t^4}-2\alpha^2\frac{1}{t^3}$  is negative for all sufficiently large $t$ and all $\alpha \neq 0$, so  the corresponding vector has critical point $\al =0$. Note that $\Ga$ does not satisfy (\ref{3.1}).
 
Needless to say, examples of such matrices abound in higher dimensions.

The following example illustrates that (\ref{3.1}) is not a sufficient condition for a vector to be associated:

 \begin{equation*} \Gamma = \frac{1}{5}\left(\begin{array}{ccc}
8&3&4\\
3&5&2\\
4&2&4\end{array}
\right);
\mbox{ then}\quad  
  \Gamma^{-1} =\frac{1}{12} \left(\begin{array}{ccc}
16&-4&-14\\

-4&16&-4\\

-14&-4&31\end{array}
\right ).\\
\end{equation*}

Where $\Gamma$ satisfies (\ref{3.1}), $\Gamma^{-1}$ is an $M$-matrix and has the first row sum negative -- hence our vector is is not associated.

 It is possible, of course, for a 3-dimensional vector to satisfy (\ref{3.1}) and have all positive row sums - i.e. be an associated vector and still have a zero entry in the inverse covariance matrix:
 \begin{equation*} \Gamma = \frac{1}{37}\left(\begin{array}{ccc}
15&4&2\\
4&6&3\\
2&3&20\end{array}
\right);
\mbox{ then}\quad  
   \Gamma^{-1} = \left(\begin{array}{ccc}
3&-2&0\\

-2&8&-1\\

0&-1&2\end{array}
\right )\\
\end{equation*}
In dimensions 4 and higher the condition (\ref{3.1}) no longer prevents the $\Ga^{-1}$ from having a zero entry at the intersection on the row sums with different signs:

\begin{equation*} \Gamma^{-1}= \left(\begin{array}{cccc}
10&-3&-3&0\\
-3&9&-2&-3\\
-3&-2&9&-3\\
0&-3&-3&\frac{65}{11} \end{array}\right),\\
\mbox{then}\quad  \Gamma= \left(\begin{array}{cccc}
\frac{257}{1400}&\frac{39}{280}&\frac{39}{280}&\frac{99}{700}\\
\frac{39}{280}&\frac{171}{616}&\frac{115}{616}&\frac{33}{140}\\
\frac{39}{280}&\frac{115}{616}&\frac{171}{616}&\frac{33}{140}\\
\frac{99}{700}&\frac{33}{140}&\frac{33}{140}&\frac{143}{350} \end{array}\right)
\end{equation*}\\

 and satisfies (\ref{3.1}).
 
 An example below is of the  3-dimensional Gaussian vector with positive, strictly positive definite covariance matrix $\Ga$ and  the diagonally dominant $\Ga^{-1}$ with one positive off diagonal entry. Hence the set of Gaussian vectors in Lemma \ref{non id} is non-vacuous already in dimension 3.
\begin{equation*} \Gamma = \frac{1}{297}\left(\begin{array}{ccc}
15&4&2\\
4&8&5\\
2&5&6\end{array}
\right);
\mbox{ then}\quad  
   \Gamma^{-1} = \left(\begin{array}{ccc}
23&-14&4\\

-14&86&-67\\

4&-67&104\end{array}
\right )\\
\end{equation*} 
  \def\noopsort#1{} \def\printfirst#1#2{#1}
\def\singleletter#1{#1}
      \def\switchargs#1#2{#2#1}
\def\bibsameauth{\leavevmode\vrule height .1ex
      depth 0pt width 2.3em\relax\,}
\makeatletter
\renewcommand{\@biblabel}[1]{\hfill#1.}\makeatother
\newcommand{\bysame}{\leavevmode\hbox to3em{\hrulefill}\,}


\newpage
\begin{thebibliography}{10}

\bibitem{B} R. Bapat, (1989) {\em Infinite divisibility of multivariate gamma distributions and M-matrices },
 Sankhya, 51, 73-78.
 
 \bibitem{F} W. Feller, {\em An introduction to Probability Theory and its Applications Vol. II},
 John Wiley and Sons, New York, 1971.





 
 \bibitem{book} M. B. Marcus and J.~Rosen, {\em Markov Processes, Gaussian
Processes and Local Times}, Cambridge studies in advanced mathematics,
100, Cambridge University Press, Cambridge, England, 2006.

\bibitem{inf div}  M. B. Marcus and J.~Rosen, {\em Infinite divisibility  of Gaussian Squares with non--zero means}, ECP, 13, (2008) 364-376.


\bibitem{critical} M. B. Marcus and J.~Rosen, {\em Existence of a critical point for the infinite  divisibility of squares  of Gaussian vectors in $R^{2}$ with non--zero mean}, Electron.J.Probab., {\em to appear}.
\bibitem{EK} N. Eisenbaum, H. Kaspi, {\em A characterization of the infinitely divisible squared Gaussian process.} Ann. Probab., 34.[6, 579], (2006).
\bibitem{GR} R.C. Griffiths,  (1984). {\em Characterizations of infinitely divisible multivariate gamma distributions. } Jour. Multivar. Anal., 15, 12-20. [6, 579]

 
 

\end{thebibliography}
\end{document}